\documentclass[10pt]{amsart}
\usepackage{amssymb}
\usepackage{bm}
\usepackage{graphicx}
\usepackage[centertags]{amsmath}
\usepackage{amsfonts}
\usepackage{amsthm}

\newtheorem{thm}{Theorem}

\newtheorem{lem}[thm]{Lemma}
\newtheorem{prop}[thm]{Proposition}

\theoremstyle{definition}

\newcommand{\rr}{\mathbb{R}}

\newcommand{\nn}{\mathbb{N}}

\newcommand{\ff}{\mathcal{F}}
\newcommand{\bb}{\mathcal{B}}

\newcommand{\dd}{\mathfrak{D}}

\newcommand{\xx}{\overline{x}}
\newcommand{\yy}{\overline{y}}
\newcommand{\zz}{\overline{z}}
\newcommand{\ww}{\overline{w}}
\newcommand{\vv}{\overline{v}}
\begin{document}
\title[A discretized approach to W. T.  Gowers' game]{A discretized approach to  W. T. Gowers' game}
\author{V. Kanellopoulos and K. Tyros}
\address{National Technical University of Athens, Faculty of Applied Sciences,
Department of Mathematics, Zografou Campus, 157 80, Athens, Greece}
\email{bkanel@math.ntua.gr, ktyros@central.ntua.gr} \footnotetext[1]{2000
\textit{Mathematics Subject Classification}: 05D10, 46B03}
\footnotetext[2]{Research supported by PEBE 2007.} \keywords{Ramsey
theory, games in Banach spaces}

\begin{abstract} We give an alternative  proof of W. T. Gowers' theorem on block
bases by reducing it to a discrete analogue on specific countable
nets. We also give a Ramsey type result on $k$-tuples of block
sequences in a normed linear space with a Schauder basis.
\end{abstract}

\maketitle


\section{Introduction}
W. T. Gowers in  \cite{G2}  (see also \cite {G1} and \cite{G3})
proved a fundamental Ramsey-type theorem for block bases in Banach
spaces which led to important discoveries in the geometry of
Banach spaces. By now there are several approaches  to Gowers'
theorem (see \cite{ADK,AT,BL1,BL2,L,R1}. Also in \cite{FFKR-N,M,
P1} there are direct proofs of  Gowers' dichotomy and in
\cite{FR,FFKR-N1,P2,R2,T} extensions and further applications).

Our aim in this note is to state and prove a discrete analogue of
Gowers' theorem which is free of approximations. To state our
results we will need the following notation. Let $\mathfrak{X}$ be a
real linear space with an infinite countable Hamel basis $(e_n)_{n}$
(actually the field over which the linear space $\mathfrak{X}$ is
defined plays no role in the arguments; it is only for the sake of
convenience that we will assume that $\mathfrak{X}$ is a real linear
space). For a subset $A\subseteq \mathfrak{X}$ by $<A>$ we denote
the linear span of $A$. Let $\mathfrak{D}$ be a subset of
$\mathfrak{X}$. By $\bb_\mathfrak{D}^\infty$ we denote the set of
all block sequences $(x_n)_n$ with $x_n\in \mathfrak{D}$ for all
$n$. For a block sequence $Z\in\bb_\mathfrak{D}^\infty$ let
$\bb_\mathfrak{D}^\infty(Z)$ be the set of all block sequences of
$\bb_\mathfrak{D}^\infty$ which are block subsequences of $Z$.

Assume that $\bb_\mathfrak{D}^\infty$ is non empty and let
$Z\in\bb_\mathfrak{D}^\infty$ and
$\mathcal{G}\subseteq\bb_\mathfrak{D}^\infty$. We define
\textit{the $\mathfrak{D}-$Gowers' game in $Z$}, denoted by
$G_\mathfrak{D}(Z)$, as follows. Player I starts the game by
choosing $W_0\in \bb_\mathfrak{D}^\infty(Z)$ and player II
responses with a vector $w_0\in <W_0>\cap \mathfrak{D}$. Then
player I chooses $W_1\in \bb_\mathfrak{D}^\infty(Z)$ and player II
chooses a vector $w_1\in <W_1>\cap \mathfrak{D}$ and so on. Player
II wins the game if the sequence $(w_0,w_1,...)$ belongs to
$\mathcal{G}$.

Suppose that  $\mathfrak{D}$ is a subset of  $ \mathfrak{X}$
satisfying the following properties.
\begin{enumerate}
\item[($\mathfrak{D}1$)] (\textit{Asymptotic property})  For all $n\in\nn$,
$\mathfrak{D}\cap <(e_i)_{i\geq n}>\neq\emptyset$.
\item[($\mathfrak{D}2$)] (\textit{Finitization property}) For all $n\in\nn$,
the set $\mathfrak{D}\cap <(e_i)_{i< n}>$ is finite.
\end{enumerate}
Property ($\mathfrak{D}1$) simply  means  that  the set of all block
sequences $\bb_\mathfrak{D}^\infty$  is non empty. Property
($\mathfrak{D}2$) implies that $\mathfrak{D}$ is countable. Hence,
endowing $\mathfrak{D}$ with the discrete topology, the space $\mathfrak{D}^\nn$
of all infinite countable sequences of $\mathfrak{D}$ equipped with the
product topology is a Polish space. We can now state our first main result.
\begin{thm} \label{main_cobinatorial_theorem}
Let $\mathfrak{X}$ be a real linear  space with a countable Hamel basis $(e_n)_n$
and let $\mathfrak{D}\subseteq \mathfrak{X}$  satisfying  properties
($\mathfrak{D}1$) and ($\mathfrak{D}2$). Also let
$\mathcal{G}\subseteq\bb_\dd^\infty$ be an analytic subset of
$\dd^\nn$. Then for every $U\in \bb_\dd^\infty$ there exists $Z\in
\bb_\dd^\infty(U)$ such that either $\bb_\dd^\infty(Z)\cap
\mathcal{G}=\emptyset$ or player II has a winning strategy in
$G_\dd(Z)$ for $\mathcal{G}$.
\end{thm}
While discrete in nature, Theorem \ref{main_cobinatorial_theorem}
can be used to derive Gowers' original result provided that
$\mathfrak{D}$ satisfies an additional property (see Section 4).

Our second main result concerns $k$-tuples of block sequences in
normed linear spaces with a Schauder basis. Precisely, let
$\mathfrak{X}$ be a real  normed linear space   with a Schauder
basis $(e_n)_n$. By $\bb_\mathfrak{X}^\infty$ we shall denote the
set of block sequences of $\mathfrak{X}$ and by
$\bb_{B_\mathfrak{X}}^\infty$ the set of all block sequences in
the unit ball $B_\mathfrak{X}$ of $\mathfrak{X}$. Two block
sequences $Z_1=(z^1_n)_n$ and $Z_2=(z^2_n)_n$ in
$\bb_\mathfrak{X}^\infty$ are said to be \textit{disjointly
supported}  if supp$z^1_n\cap $supp$z^2_m=\emptyset$ for all
$m,n$. For a positive integer $k\geq 2$ and for every
$Z\in\bb_\mathfrak{X}^\infty$, the set of all $k$-tuples
consisting of pairwise disjointly supported block subsequences of
$Z$ in $B_\mathfrak{X}$ will be denoted by
$(\bb_{B_\mathfrak{X}}^\infty(Z))_\perp^k$. Also, for a family
$\mathfrak{F}\subseteq (\bb_\mathfrak{X}^\infty)^k$ of $k$-tuples
of block sequences of $\mathfrak{X}$, the \textit{upwards closure}
of $\mathfrak{F}$ is defined to be  the set
\begin{eqnarray*}
\mathfrak{F}^\uparrow=\big\{(U_i)_{i=0}^{k-1}\in (\bb_\mathfrak{X}^\infty)^k
& : & \exists  (V_i)_{i=0}^{k-1}\in\mathfrak{F} \text{ such that}\\
& & \forall i \; V_i \text{ is a block subsequence of } U_i\big\}
\end{eqnarray*}
If $\Delta=(\delta_n)_n$ is a sequence of positive reals, then
\textit{the $\Delta$-expansion of $\mathfrak{F}$} is defined to be
the set
\[\mathfrak{F}_{\Delta}=\big\{(U_i)_{i=0}^{k-1}\in
(\bb_\mathfrak{X}^\infty)^k:\exists
(V_i)_{i=0}^{k-1}\in\mathfrak{F}\text{ such that}
\;\forall i\; dist(U_i,V_i)\leq\Delta\big\}.\]
We prove the following.
\begin{thm}\label{Ramsey}
Let $\mathfrak{X}$ be a real normed linear space with a Schauder
basis, $k\geq 2$ and $\mathfrak{F}$ be an analytic subset of
$(\bb_{B_\mathfrak{X}}^\infty)^k$. Then for every sequence of
positive real numbers $\Delta=(\delta_n)_n$ there is
$Y\in\bb_{\mathfrak{X}}^\infty$ such that either
$(\bb_{B_\mathfrak{X}}^\infty(Y))_\perp^k\cap\mathfrak{F}=\emptyset$
or $(\bb_{B_\mathfrak{X}}^\infty(Y))^k\subseteq
(\mathfrak{F}_{\Delta})^\uparrow$.
\end{thm}
In the above theorem the topology of $\bb_{B_\mathfrak{X}}^\infty$
is the induced one by the product of the norm topology. Theorem
\ref{Ramsey} applied for k=2 and the family
\[ \mathfrak{F}=\{(U_1,U_2)\in (\bb_{B_\mathfrak{X}}^\infty)^2: U_1,U_2
\;\text{are $C$- equivalent}\}\] where $C\geq 1$ is a constant, yields Gowers' second
dichotomy (see Lemma 7.3 in \cite{G2}).

\section{Notation.} Let $\mathfrak{X}$ be a real linear
space with an infinite countable Hamel basis $(e_n)_n$. For two
non zero vectors $x,y$ in $\mathfrak{X}$, we write $x<y$ if $\max$
supp $x<\min$ supp $y$, (where supp $x$ is the \textit{support} of
$x$, i.e. if $x=\sum_n\lambda_n e_n$ then supp
$x=\{n\in\nn:\lambda_n\neq 0\}$). A sequence $(x_n)_n$ of vectors
in $\mathfrak{X}$ is called a \textit{block sequence} (or
\textit{block basis}) if $x_n<x_{n+1}$ for all $n$.

 Capital letters (such us $U,V,Y, Z, ...$) refer to
infinite block sequences and lower case  letters with a line over
them (such us $\overline{u},\vv,\yy, \zz, ...$) to finite block
sequences. We write $Y\preceq Z$ to denote that $Y$ is a
\textit{block subsequence} of $Z$,  that is $Y=(y_n)_n$, $Z=(z_n)_n$
are block sequences  and for all $n$, $y_n\in<(z_i)_i>$. The
notation $\yy\preceq Z$ and $\yy\preceq\zz$ are defined analogously.
For $\xx=(x_n)_{n=0}^k$ and $Y=(y_n)_n$  we write $\xx<Y$, if
$x_k<y_0$. For $\xx<Y$, $\xx^\smallfrown Y$ denotes the block
sequence $(z_n)_n$ that starts with the elements of $\xx$ and
continues with these of $Y$. Also for $\xx<\yy$,   the finite block
sequence $\xx^\smallfrown\yy$ is similarly defined. For a block
sequence $Z=(z_n)_n$ and an infinite subset $L$ of $\nn$  we set
$Z|_L=(z_n)_{n\in L}$. Also for $k\in\nn$, $Z|_k=(z_n)_{n=0}^{k-1}$
(where for $k=0$, $Z|_0=\emptyset$).

Let $\mathfrak{D}$ be a subset of  $ \mathfrak{X}$. By
$\bb_\mathfrak{D}^\infty$ (resp. $\bb_\mathfrak{D}^{<\infty}$) we
denote the set of all infinite  (resp. finite) block sequences
$(x_n)_n$ with $x_n\in \mathfrak{D}$ for all $n$. The set of all
infinite (resp. finite) block sequences in $\mathfrak{X}$ is denoted
by  $\bb^\infty_{\mathfrak{X}}$ (resp.
$\bb^{<\infty}_\mathfrak{X}$). For $Z\in\bb^\infty_\mathfrak{X}$ we
set $\bb_\mathfrak{D}^\infty(Z)=\{Y\in
\bb_\mathfrak{D}^\infty:Y\preceq Z\}$ and
$\bb_\mathfrak{D}^{<\infty}(Z)=\{\yy\in\bb_\mathfrak{D}^{<\infty}:\yy\preceq
Z\}$. Similarly for $\zz\in \bb_\mathfrak{X}^{<\infty}$,
$\bb_\mathfrak{D}^{<\infty}(\zz)=\{\yy\in\bb_\mathfrak{D}^{<\infty}:\yy\preceq
\zz\}$. For a block sequence $Z\in\bb_\mathfrak{D}^{\infty}$, we set
$< Z >_{\mathfrak{D}}=< Z >\cap \mathfrak{D}$ where $<Z>$ is the
linear span of $Z$.

\section{Discretization of Gowers' game.}
Throughout this section,  $\mathfrak{X}$ is a real linear space with
countable Hamel basis $(e_n)_n$ and $\mathfrak{D}$  is a subset of
$\mathfrak{X}$ satisfying properties ($\mathfrak{D}1$) and
($\mathfrak{D}2$) as stated in the Introduction. Notice that
($\mathfrak{D}2$)  also gives that for every $U=(u_i)_i\in
\bb_\dd^\infty$ and  $n\in\nn$, the set
$\bb_\dd^{<\infty}((u_i)_{i<n})$ is finite.
\subsection{Admissible families of $\mathfrak{D}$-pairs.} The aim of this subsection is to review the  methods
that we will follow to handle  the  several diagonalizations that
will appear (see also \cite{G2}, \cite{PR}).
 A $\mathfrak{D}$-\textit{pair}  is a pair
$(\xx,Y)$ where $\xx\in\bb_\mathfrak{D}^{<\infty}$ and
$Y\in\bb_\mathfrak{D}^\infty$.  For  $U\in\bb_\dd^\infty$, a family
$\mathcal{P}\subseteq \bb^{<\infty}_\dd(U)\times \bb^\infty_\dd(U)$
is called \textit{admissible family of $\mathfrak{D}$- pairs in} $U$
if it satisfies the next properties:
\begin{enumerate}
\item[($\mathcal{P}1$)] (\textit{Heredity})  If $(\xx,Y)\in
\mathcal{P}$ and $Z\in\bb_\dd^\infty(Y)$ then $(\xx,
Z)\in\mathcal{P}$.
\item[($\mathcal{P}2$)] (\textit{Cofinality}) For every $(\xx, Y)\in
\bb^{<\infty}_\dd(U)\times \bb^\infty_\dd(U)$, there is
$Z\in\bb_\dd^\infty(Y)$ such that $(\xx,Z)\in \mathcal{P}$.
\end{enumerate}
For simplicity in  the sequel  when we write  ``pair" we will always
mean a ``$\mathfrak{D}$-pair". It will  often happen that an
admissible family of pairs  has one more property.
\begin{enumerate}
\item[($\mathcal{P}3$)]  If $(\xx,Y)\in \mathcal{P}$, $\xx< Y$ and  $k=\min\{m:\xx\in\bb_\mathfrak{D}^{<\infty}((u_i)_{i=1}^m)\}$
 then for every $\yy\in
 \bb_\mathfrak{D}^{<\infty}((u_i)_{i=1}^k)$,
 $(\xx,\yy^\smallfrown Y)\in\mathcal{P}.$
\end{enumerate}
The next  lemma follows by a standard diagonalization
 argument.

\begin{lem}
  \label{diagonalization 2'}
Let $U\in\bb_\dd^\infty$ and let $\mathcal{P}$ be an admissible
family of pairs in $U$. Then there is
$W\in\bb_\mathfrak{D}^\infty(U)$ such that for all
$\ww\in\bb_\mathfrak{D}^{<\infty}(W)$ and all
$Y\in\bb_\dd^\infty(W)$ with $\ww< Y$, $(\ww,Y)\in\mathcal{P}.$ If
in addition  $\mathcal{P}$  satisfies ($\mathcal{P}3$)  then for all
$\ww\in\bb_\mathfrak{D}^{<\infty}(W)$, $(\ww,W)\in\mathcal{P}.$
\end{lem}




\subsection{The discrete Gowers' game.}

Given  $Y\in\bb_\mathfrak{D}^\infty$ and a family of infinite
block sequences $\mathcal{G}\subseteq\bb_\mathfrak{D}^\infty$, we
define the \textit{$\mathfrak{D}-$Gowers' game},
$G_\mathfrak{D}(Y)$, as follows. Player I starts the game by
choosing $Z_0\in \bb_\mathfrak{D}^\infty(Y)$ and player II
responses with a vector $z_0\in <Z_0>_\mathfrak{D}$. Then player I
chooses $Z_1\in \bb_\mathfrak{D}^\infty(Y)$ and player II chooses
a vector $z_1\in <Z_1>_\mathfrak{D}$ with $z_0<z_1$ and so on.
More generally for a finite block sequence $\xx\in
\bb_\mathfrak{D}^{<\infty}$ and $Y\in\bb_\mathfrak{D}^\infty$ the
game $G_\mathfrak{D}(\xx, Y)$ is defined as above with the
additional condition that player II in the first move chooses
$z_0>\xx$. Clearly $G_\dd(\emptyset, Y)$ is identical to
$G_\dd(Y)$. We will say that player II \textit{wins the game}
$G_\dd(\xx,Y)$ \textit{for} $\mathcal{G}$ if the block sequence
$\xx^\smallfrown(z_0,z_1,...)$ belongs to $\mathcal{G}$.

The basic terminology that we shall use is an adaptation of the
classical Galvin- Prikry's one (cf. \cite{GP}, \cite{E}) in the frame of Gowers' game. More
precisely, for $\xx\in\bb_\mathfrak{D}^{<\infty}$,
$Y\in\bb_\mathfrak{D}^\infty$ and
$\mathcal{G}\subset\bb_\mathfrak{D}^\infty$ we say that $Y$
$\mathcal{G}-$\textit{ accepts} $\xx$ if player II has a winning
strategy in $G_\mathfrak{D}(\xx, Y)$ for $\mathcal{G}$ and that $Y$
$\mathcal{G}-$\textit{ rejects} $\xx$ if there is no
$Z\in\bb_\mathfrak{D}^\infty(Y)$ which $\mathcal{G}-$ accepts $\xx$.
We also say that $Y$ $\mathcal{G}-$\textit{ decides} $\xx$ if either
$Y$ $\mathcal{G}$- accepts $\xx$ or $Y$ $\mathcal{G} $-rejects
$\xx$.

Notice that
if $\xx=\emptyset$ then to say that ``$Y$ $\mathcal{G}$-accepts the
empty sequence" means that  player II has a winning strategy in
$G_\mathfrak{D}(Y)$ for $\mathcal{G}$. Similarly the statement that
``$Y$ $\mathcal{G}$-rejects the empty sequence" is equivalent to
that for all $Z\in\bb_\mathfrak{D}^{\infty}(Y)$ player II has no
winning strategy in $G_\mathfrak{D}(Z)$ for $\mathcal{G}$. The
following lemma is  easily verified.
\begin{lem} \label{decide}For every $U\in\bb_\mathfrak{D}^\infty$
 and every $\mathcal{G}\subseteq \bb_\dd^\infty$,
the family \[\mathcal{P}=\{(\xx,Y)\in
\bb_\mathfrak{D}^{<\infty}(U)\times\bb^\infty_\mathfrak{D}(U):\;Y
\;\mathcal{G}-\text{decides}\;\xx\}\] is an admissible family of
pairs in $U$ which in addition satisfies property  ($\mathcal{P}3$).
\end{lem}
Actually the family $\mathcal{P}$ of the above lemma satisfies the
following  stronger than ($\mathcal{P}3$)  property: If
$(\xx,Y)\in\mathcal{P}$ and  $Z\in\bb_\mathfrak{D}^{\infty}$ such
that  there is $n\in\nn$ with $Z|_{[n,\infty)}\preceq Y$, then
$(\xx, Z)\in \mathcal{P}$.

 For the sake of simplicity in the following we will
omit the letter $\mathcal{G}$ in front of the words ``accepts",
``rejects" and ``decides". The next lemma is a consequence of Lemma
\ref{decide} and Lemma \ref {diagonalization 2'}.
\begin{lem}
  \label{first_combinatorial_lemma_for_the_open}
For every  $U\in \bb_\mathfrak{D}^\infty$ there is
$W\in\bb_\mathfrak{D}^\infty(U)$  such that for all
$\ww\in\bb_\mathfrak{D}^{<\infty}(W)$,
 $W$ decides $\ww$.
\end{lem}
The crucial point at which the above notions of ``accept-reject"
essentially differ from the original ones reveals in the next lemma.
Here the notion of the winning strategy  replaces successfully the
traditional pigeonhole principle.
\begin{lem}
  \label{second_combinatorial_lemma_for_the_open}
  Let $W\in\bb_\mathfrak{D}^\infty$ such that $W$ decides  all
$\ww\in\bb_\mathfrak{D}^{<\infty}(W)$ and assume that there is
$\ww_0\in\bb_\mathfrak{D}^\infty(W)$ such that $W$ rejects $\ww_0$.
Then for every $Y\in\bb_\mathfrak{D}^\infty(W)$ there is
$Z\in\bb_\mathfrak{D}^\infty(Y)$ such that for every
$z\in<Z>_\mathfrak{D}$ with $\ww_0<z$, $W$ rejects $\ww_0^\frown z$.
\end{lem}
\begin{proof}
If the conclusion is false then there is
$Y\in\bb_\mathfrak{D}^\infty(W)$ such that for every
$Z\in\bb_\mathfrak{D}^\infty(Y)$ there is  $z\in<Z>_\mathfrak{D}$
with $\ww_0<z$  such that  $W$ accepts $\ww_0^\frown z$. It is
easy to see that this means that  player II has a winning strategy
in $G_\mathfrak{D}(\ww_0, Y)$ for $\mathcal{G}$ and thus $Y$ accepts
$\ww_0$. But this is a contradiction since
$Y\in\bb_\mathfrak{D}^\infty(W)$ and $W$ rejects $\ww_0$.
\end{proof}
\begin{lem} \label{combinatorial_proposition_for_open_family} For every
$U\in\bb_\mathfrak{D}^\infty$  there exists
$Z\in\bb_\mathfrak{D}^\infty(U)$ such that either $Z$ rejects all
$\zz\in \bb_\mathfrak{D}^{<\infty}(Z)$ or player II has winning
strategy in $G_\mathfrak{D}(Z)$ for $\mathcal{G}.$
\end{lem}
\begin{proof}
 By Lemma
\ref{first_combinatorial_lemma_for_the_open} there is
  $W\in \bb_\mathfrak{D}^\infty(U)$ such that for every $\ww\in \bb_\mathfrak{D}^{<\infty}(W)$, $W$ decides $\ww$.
If $W$ accepts the empty sequence then we readily have the second
alternative of the conclusion for $Z=W$. In the opposite case consider the
following family  in $ \bb^{<\infty}_\dd(W)\times
\bb^\infty_\dd(W)$:
\[\mathcal{P}=\{(\xx, Y):\text{Either }
W\;\text{accepts}\; \xx\;\text{or} \;\forall\;
y\in<Y>_\dd\;\text{with} \;\xx<y, \;W\;\text{rejects
} \xx^\smallfrown y\}\] Using Lemma \ref{second_combinatorial_lemma_for_the_open} we
easily verify that $\mathcal{P}$ is an admissible family in $W$
which satisfies also property ($\mathcal{P}3$). Hence by Lemma
\ref{diagonalization 2'} there is $Z\in\bb_\dd^\infty(W)$ such that
for every $\zz\in\bb^{<\infty}_\dd(Z)$, $(\zz,Z)\in\mathcal{P}$. By
our assumption $W$ rejects the empty sequence. Hence since
$(\emptyset,Z)\in\mathcal{P}$ we have that $W$ and so $Z$ rejects
all $z\in <Z>_\dd$. By induction on the length of finite block
sequences in $\bb_\dd^{<\infty}(Z)$, it is easily shown that  $Z$
rejects all $\zz\in \bb_\dd^{<\infty}(Z)$. \end{proof}

We have finally arrived at our first stop which is an analog of the
well known result of Nash-Williams (\cite{N-W}). Consider the set $\mathfrak{D}$
as a topological space with the discrete topology and
$\mathfrak{D}^\nn$ with the product topology.
\begin{lem}\label{lemma 4}
Let  $\mathcal{G}\subseteq \bb_\mathfrak{D}^\infty$ be  open in  $\mathfrak{D}^\nn$. Then for every
$U\in\bb_\mathfrak{D}^\infty$ there exists
$Z\in\bb_\mathfrak{D}^\infty(U)$ such that either
$\bb_\mathfrak{D}^\infty(Z)\cap\mathcal{G}=\emptyset$ or
  player II has  a winning strategy in $G_\mathfrak{D}(Z)$ for $\mathcal{G}$.
\end{lem}
\begin{proof}
By Lemma  \ref{combinatorial_proposition_for_open_family} we can
find $Z\in\bb_\mathfrak{D}^\infty(U)$ such that either $Z$ rejects
all $\zz\in\bb_\mathfrak{D}^{<\infty}(Z)$, or player II has a
winning strategy in $G_\mathfrak{D}(Z)$ for $\mathcal{G}$. Hence it
suffices to show that   the first alternative gives that
$\bb_\mathfrak{D}^\infty(Z)\cap\mathcal{G}=\emptyset$. Indeed, let
$W=(w_n)_{n}\in\bb_\mathfrak{D}^\infty(Z)$. Then for all  $k$, $Z$
rejects $W|_k=(w_n)_{n<k}$. Therefore
there is some $Z_k\in\bb_\mathfrak{D}^\infty(Z)$ with $W|_k<Z_k$
such that  $W|_k^\smallfrown Z_k\not\in\mathcal{G}$. Since the
sequence $(W|_k^\smallfrown Z_k)_{k}$ converges in
$\mathfrak{D}^\nn$ to $W$ and the complement of $\mathcal{G}$ is
closed, we conclude  that $W\not\in\mathcal{G}$.
\end{proof}

We pass now to the  case of an analytic family $\mathcal{G}$.
First let us state  some basic definitions (cf. \cite{K}). Let  $\nn^{<\nn}$ be the set of all finite sequences in $\nn$ and
let $\mathcal{N}$ be  the Baire space i.e. the space of all infinite sequences in $\nn$ with the topology generated by the sets
$\mathcal{N}_s=\{\sigma\in \mathcal{N}:\;\exists
n\;\text{with}\;\sigma|n=s\}$, $s\in\nn^{<\nn}$. A subset  of a Polish space $X$ is called \textit{analytic}
 if it is the image of a continuous function from  $\mathcal{N}$ into $X$.

For the
next  lemmas we fix the following.
\begin{enumerate}
\item [(a)] A family $(\mathcal{G}^s)_{s\in\nn^{<\nn}}$ of subsets of
$\bb_\dd^\infty$ such that for all $s$,
$\mathcal{G}^s=\bigcup_n\mathcal{G}^{s^\smallfrown n}$.
\item[(b)]
A bijection $\varphi:\nn^{<\nn}\to \nn$ such that
$\varphi(\emptyset)=0$ and
for all $s,n$, $\varphi(s^\smallfrown n)>\varphi(s)$.
\end{enumerate}
For each $\xx$ in $\bb_\mathfrak{D}^{<\infty}$ we set $s_{\xx}$ to
be the unique element element of $\nn^{<\nn}$ such that
$\varphi(s_{\xx})$ equals to the length of $\xx$.
 For a
$\mathfrak{D}$- pair $(\xx,Y)$ we set
\[\bb_\mathfrak{D}^\infty(\xx,Y)=\{V\in\bb_\mathfrak{D}^\infty:
\exists k\;\text{such that}\;
V|_k=\xx\;\text{and}\;V|_{[k,\infty)}\preceq Y\}\] Finally, recall
the following terminology from \cite{G2}. For a family
$\mathcal{G}\subseteq \bb^\infty_\mathfrak{D}$ we say that
$\mathcal{G}$ is \textit{large for} $(\xx, Y)$ if for all
$Z\in\bb^\infty_\mathfrak{D}(Y)$, $\mathcal{G}\cap
\bb^\infty_\mathfrak{D}(\xx,Z)\neq\emptyset$. In the case
$\xx=\emptyset$ we simply say that $\mathcal{G}$ is large for $Y$.

\begin{lem}
\label{first_combinatorial_lemma_for_the_analytic_family} For
every $U\in\bb^\infty_\mathfrak{D}$  there is
$W\in\bb_\mathfrak{D}^\infty(U)$ such that for every $ \ww \in
\bb_\mathfrak{D}^{<\infty}(W)$, either $\mathcal{G}^{s_{\ww}}\cap
\bb^\infty_\mathfrak{D}(\ww,W)=\emptyset$ or
$\mathcal{G}^{s_{\ww}}$ is large for $(\ww,W)$.
\end{lem}

\begin{proof}
Let $\mathcal{P}$ be the set of all pairs $(\xx,Y)$ in
$\bb_{\mathfrak{D}}^{<\infty}(U)\times\bb_{\mathfrak{D}}^\infty(Y)$
such that either $\mathcal{G}^{s_{\xx}}\cap
\bb^\infty_\mathfrak{D}(\xx,Y)=\emptyset$ or
$\mathcal{G}^{s_{\xx}}$ is large for $(\xx,Y)$. It is easy to see
that $\mathcal{P}$ is admissible satisfying property
($\mathcal{P}3)$. Hence the conclusion follows by Lemma
\ref{diagonalization 2'}.
\end{proof}

Let $W\in\bb_\dd^\infty$ be  a block sequence in $\mathfrak{D}$
 satisfying the conclusion of Lemma
\ref{first_combinatorial_lemma_for_the_analytic_family}. For
$\ww\in\bb_\mathfrak{D}^{<\infty}(W)$, let $\mathcal{F}(\ww)$ be
the family of all $V=(v_i)_i\in\bb^{\infty}_\mathfrak{D}( W)$ with
$\ww<V$ and the following properties. There exist $m, l\in\nn$
with $l\geq 1$ such that
\begin{enumerate}
\item[(i)] $s_{\ww}^\smallfrown m=s_{\xx}$, where $\xx=
\ww^\smallfrown(v_i)_{i=0}^{l-1}$ and \item[(ii)] The family
$\mathcal{G}^{s_{\ww}^\smallfrown m}$ is large for
$(\ww^\smallfrown (v_i)_{i=0}^{l-1}, W)$.
\end{enumerate}
Notice that $\mathcal{F}( \ww)$ is open in $\mathfrak{D}^\nn$.
\begin{lem}
  \label{second__combinatorial_lemma_for_the_analytic_family}
Let  $\ww\in \bb^{<\infty}_\mathfrak{D}(W)$ and assume that
$\mathcal{G}^{s_{\ww}}$ is large for $(\ww,W)$. Then
$\mathcal{F}(\ww)$ is  large for $W$.
\end{lem}
\begin{proof}
Let $Z\in\bb_\mathfrak{D}^\infty(W)$.  Since
$\mathcal{G}^{s_{\ww}}$ is large for $(\ww,W) $ there is
$V=(v_i)_i$ such that $\ww<V$ and $\ww^\smallfrown V\in
\mathcal{G}^{s_{\ww}}\cap \bb_\mathfrak{D}^\infty(\ww,Z)=\bigcup_m
\mathcal{G}^{s_{\ww}^\smallfrown
m}\cap\bb_\mathfrak{D}^\infty(\ww,Z)$ and so for some $m\in\nn$,
$\ww^\smallfrown V\in \mathcal{G}^{s_{\ww}^\smallfrown
m}\cap\bb_\mathfrak{D}^\infty(\ww,Z)$. Notice that for
$l=\varphi(s^\smallfrown m)-\varphi(s)$ we have that
$s_{\ww}^\smallfrown m=s_{\xx}$, where $\xx=\ww^\smallfrown
(v_i)_{i=0}^{l-1}$, and  $\ww^\smallfrown V\in
\mathcal{G}^{s_{\ww}^\smallfrown m}\cap
\bb^\infty_\mathfrak{D}(\ww^\smallfrown (v_i)_{i=0}^{l-1},Z)$.
Therefore  $\mathcal{G}^{s_{\ww}^\smallfrown m}\cap
\bb^\infty_\mathfrak{D}(\ww^\smallfrown
(v_i)_{i=0}^{l-1},W)\neq\emptyset$ which (by the properties of
$W$) means that $\mathcal{G}^{s_{\ww}^\smallfrown m}$ is large for
$(\ww^\smallfrown (v_i)_{i=0}^{l-1},W)$. Hence
$V\in\mathcal{F}(\ww)\cap\bb_\mathfrak{D}^\infty(Z)$.\end{proof}

\begin{lem}\label{third__combinatorial_lemma_for_the_analytic_family}
There is $Z\in\bb_\dd^\infty(W)$ such that for every $\zz\in
 \bb^{<\infty}_\mathfrak{D}(Z)$ we have that either
$\mathcal{G}^{s_{\ww}}\cap \bb_\dd^\infty(\zz,Z)=\emptyset$ or
player II has a winning strategy in the game $G_\dd(Z)$ for the
family $\mathcal{F}(\zz)$.
\end{lem}
\begin{proof}
 Let  $\mathcal{P}$ be the family
 of pairs $(\ww, Y)\in
\bb_\dd^{<\infty}(W)\times \bb_\dd^\infty(W)$  such that either
$\mathcal{G}^{s_{\ww}}\cap\bb^\infty_\mathfrak{D}(\ww,Y)=\emptyset$
or
 player II has a
winning strategy in the game $G_\dd(Y)$ for the family
$\mathcal{F}(\ww)$.

By Lemma \ref{diagonalization 2'} it suffices to show  that
$\mathcal{P}$ is an admissible family of pairs in $W$ which in
addition satisfies property $(\mathcal{P}3)$. It is easy to see
that only the cofinality property needs some explanation. To this
end let $(\ww,Y)\in\bb^{<\infty}_\dd(W)\times\bb^{\infty}_\dd(W)$.
Since $\ww\in \bb_\dd^{<\infty}(W)$ we have that either
$\mathcal{G}^{s_{\ww}}\cap \bb_\dd^\infty(\ww,W)=\emptyset$, or
$\mathcal{G}^{s_{\ww}}$ is large  for $(\ww,W)$. In the first
case, $\mathcal{G}^{s_{\ww}}\cap \bb_\dd^\infty(\ww,Y)=\emptyset$
and so $(\ww,Y)\in\mathcal{P}$. In the second case,  Lemma
\ref{second__combinatorial_lemma_for_the_analytic_family} implies
that $\mathcal{F}(\ww)$ is  large for $W$. Hence by Lemma
\ref{lemma 4}, there is $V\in\bb_\dd^\infty(Y)$ such that player
II has a winning strategy in $G_\mathfrak{D}(V)$ for
$\mathcal{F}(\ww)$ and so $(\ww,V)\in\mathcal{P}$.\end{proof}

We are now ready for the proof of the main  result.
\medskip

\begin{proof} [\textbf{Proof of Theorem \ref{main_cobinatorial_theorem}:
}]
Assume that there is no $Z\in\bb_\dd^\infty(U)$ such that
$\bb_\dd^\infty(Z)\cap \mathcal{G}=\emptyset$, that is
$\mathcal{G}$ is large for $U$. Let $f:\mathcal{N}\to
\mathfrak{D}^{\nn}$ be a continuous map with
$f[\mathcal{N}]=\mathcal{G}$ and for $s\in\nn^{<\nn}$, let
$\mathcal{G}^s=f[\mathcal{N}_s]$. Then
$\mathcal{G}^\emptyset=\mathcal{G}$ and
 $\mathcal{G}^s=\bigcup_n\mathcal{G}^{s^\smallfrown n}$.
Following the process of  the above lemmas let
$W\in\bb_\dd^\infty(U)$  be  as in   Lemma
\ref{first_combinatorial_lemma_for_the_analytic_family} and  $Z\in\bb_\dd^\infty(W)$ as in   Lemma
\ref{third__combinatorial_lemma_for_the_analytic_family}. We claim
that player II has a winning strategy in the game $G_\dd(Z)$ for $\mathcal{G}$.

Indeed,  by our assumption  $\mathcal{G}=\mathcal{G}^\emptyset$ is
large in $\bb_\dd^\infty(Z)=\bb_\dd^\infty(\emptyset,Z)$ and so
player II has a winning strategy in $G_\dd(Z)$ for
$\mathcal{F}(\emptyset)$. This means that player II is able  to
produce after a finite number of moves, a finite block sequence
$\yy_0\in\bb_\dd^{<\infty}(Z)$ such that there is $m_0\in\nn$,
with  $s_{\yy_0}=(m_0)$ and $\mathcal{G}^{( m_0)}$ large for
$(\yy_0,W)$. By Lemma
\ref{third__combinatorial_lemma_for_the_analytic_family},
 player II has a winning
strategy   in $G_\dd(Z)$ for $\mathcal{F}(\yy_0)$, that is player
II can extend $\yy_0$ to a finite block sequence
$\yy_0^\smallfrown \yy_1\in\bb_\dd^{<\infty}(Z)$ such that there
is $m_1\in\nn$ such that $ s_{\yy_0^\smallfrown \yy_1}=(m_0, m_1)$
 and $\mathcal{G}^{(m_0, m_1)}$ is large for
$(\yy_0^\smallfrown \yy_1,W)$.

Continuing in this way we conclude that player II has a strategy
in the game $G_\dd(Z)$ to construct a block sequence
$Y=\yy_0^\smallfrown\yy_1^\smallfrown ...$ such that for some
$\sigma=(m_i)_i\in\mathcal{N}$ and
 for every $k\in\nn$, $\mathcal{G}^{\sigma|k}$ is large for
$((\yy_0^\smallfrown...^\smallfrown \yy_{k-1}),W)$. To show that
 this is actually  a winning strategy for $\mathcal{G}$ we have to
prove that $Y\in\mathcal{G}$. Fix $k\in\nn$. Since
$\mathcal{G}^{\sigma|k}$ is large for
$((\yy_0^\smallfrown...^\smallfrown \yy_{k-1}),W)$, we have that there
exists  $Y_k\in\bb_\dd^\infty(W)$ such that
$(\yy_0^\smallfrown...^\smallfrown \yy_{k-1})^\smallfrown Y_k\in
\mathcal{G}^{\sigma|k}$. Since $(\mathcal{G}^{\sigma|n})_n$ is decreasing,
 $Y=\lim_n (\yy_0^\smallfrown....^\smallfrown \yy_{n-1})^\smallfrown
Y_n\in \overline{\mathcal{G}^{\sigma|k}}$, for all $k\in\nn$, and
thus $Y\in \cap_k\overline{\mathcal{G}^{\sigma|k}}$. By the
continuity of $f$, $\cap_k\overline{\mathcal{G}^{\sigma|k}}
=\{f(\sigma)\}$ and therefore $Y=f(\sigma)\in\mathcal{G}$.
\end{proof}
  \section{Passing from the discrete  to Gowers' game.}
In this section we will see how using Theorem
\ref{main_cobinatorial_theorem} one can derive W. T. Gowers'
Ramsey theorem (see Theorem
\ref{main_Gowers_combinatorial_theorem}). From now on and for all
the rest of this note $\mathfrak{X}$ will be a normed linear space
with a Schauder basis $(e_n)_n$.

First let us recall some relevant definitions. Let
$\bb_\mathfrak{X}^\infty$ (resp. $\bb_{B_\mathfrak{X}}^\infty$) be
the set of all block sequences in $\mathfrak{X}$ (resp. in the
unit ball $B_\mathfrak{X}$ of $\mathfrak{X}$ ). Let
$U=(u_n)_{n},V=(v_n)_{n}\in \bb_\mathfrak{X}^\infty$ and
$\Delta=(\delta_n)_{n}$ a sequence of positive real numbers. We
say that $U,V$ are $\Delta-$\textit{near} and we write
$dist(U,V)\leq\Delta$ if for all $n\in\nn$,
$\|u_n-v_n\|\leq\delta_n$. For a family
$\ff\subseteq\bb_\mathfrak{X}^\infty$  and
 a sequence $\Delta=(\delta_n)_{n}$ of positive real numbers the
$\Delta$-\textit{expansion} of $\mathcal{F}$ is the set
\[\ff_\Delta=\{U\in\bb_{\mathfrak{X}}^\infty:\exists V\in\ff\text{ such that
}dist(U,V)\leq\Delta\}\] For $Y\in\bb^\infty_{B_\mathfrak{X}}$ and
a family $\ff\subseteq \bb^\infty_{B_\mathfrak{X}}$ the Gowers'
game $G_{\mathfrak{X}}(Y)$ is defined  as the
$\mathfrak{D}$-Gowers game by replacing $\mathfrak{D}$ and
$\mathcal{G}\subseteq \bb_\mathfrak{D}^\infty $ with the unit ball
$B_\mathfrak{X}$  and  $\ff\subseteq \bb_{B_\mathfrak{X}}^\infty$
respectively.

 For the next two
lemmas we fix the following. \begin{enumerate}\item[(i)] A subset
$\mathfrak{D}$ of $<(e_n)_n>$ satisfying
 the asymptotic property ($\mathfrak{D}1$).
 \item[(ii)] A family $\ff\subseteq \bb_{B_\mathfrak{X}}^\infty$  of block sequences
in $B_\mathfrak{X}$,   \item[(iii)] A sequence
$\Delta=(\delta_n)_{n}$ of positive real numbers.\end{enumerate}

\begin{lem}\label{lemma2_from_discrete_to_continuous} Let
$\mathcal{G}=\ff_{\Delta}\cap\bb_\mathfrak{D}^\infty$ and suppose that for some $\widetilde{Z}\in\bb_\dd^\infty$,
 $\bb_\mathfrak{D}^\infty(\widetilde{Z})\cap \mathcal{G}=\emptyset$. Assume that
there exist $Z\in\bb_\mathfrak{X}^\infty$  such that
\[\bb_{B_\mathfrak{X}}^\infty
(Z)\subseteq(\bb_\mathfrak{D}^\infty(\widetilde{Z}))_{\Delta}\]
(that is for every  block subsequence $U=(u_n)_n$ of $Z$ with
$\|u_n\|\leq 1$ there is a block subsequence
$\widetilde{U}=(\widetilde{u}_n)_n$ of $\widetilde{Z}$ with
$\widetilde{u}_n\in \mathfrak{D}$ such that
$dist(U,\widetilde{U})\leq\Delta$).

Then $\bb_{B_\mathfrak{X}}^\infty(Z)\cap\ff=\emptyset$.\end{lem}
\begin{proof}
Let $U\in\bb_{B_\mathfrak{X}}^\infty(Z)$. By our assumptions there
is $\widetilde{U}\in\bb_\mathfrak{D}^\infty(\widetilde{Z})$ such
that $dist(U,\widetilde{U})\leq\Delta$ and
$\widetilde{U}\not\in\mathcal{G}$. Then $U\not\in\ff$, otherwise
$\widetilde{U}\in \ff_\Delta\cap\bb_\mathfrak{D}^\infty(\widetilde{Z})$
which is a contradiction.
  \end{proof}
\begin{lem}\label{lemma1_from_discrete_to_continuous}
Let   $\delta_0\leq 1$ and $\ \sum_{j=n+1}^\infty\delta_j\leq
\delta_n,$ for all $n$.  Let
$\mathcal{G}=\ff_{\Delta/10C}\cap\bb_\mathfrak{D}^\infty$, where $C$
is the basis constant of $(e_n)_n$ and suppose that for some $\widetilde{Z}\in\bb_\dd^\infty$
 player
II has a winning strategy in the discrete  game
$G_\mathfrak{D}(\widetilde{Z})$ for $\mathcal{G}$.
Assume that there exist
$Z\in\bb_\mathfrak{X}^\infty$
 such that
\[\bb_{B_\mathfrak{X}}^\infty
(Z)\subseteq(\bb_\mathfrak{D}^\infty(\widetilde{Z}))_{\Delta/10 C}\]
Then  player
II has a winning strategy in Gowers' game $G_\mathfrak{X}(Z)$ for $\ff_\Delta$.
\end{lem}
\begin{proof}
We will define  a winning strategy for player II in
Gowers' game $G_\mathfrak{X}(Z)$ for $\ff_\Delta$ provided that he
has one in the discrete game $G_\mathfrak{D}(Z)$ for $\mathcal{G}$.
Suppose that we have just completed  the n-th move of the game
$G_\mathfrak{X}(Z)$ (resp. of the discrete game $G_\mathfrak{D}(\widetilde{Z})$) and  $x_0<...<x_{n-1}$  (resp.
$\widetilde{x}_0<...<\widetilde{x}_{n-1}$) have been chosen by
player II in $G_\mathfrak{X}(Z)$ (resp. in
$G_\mathfrak{D}(\widetilde{Z})$).

Suppose that in the game $G_\mathfrak{X}(Z)$ player  I chooses a
block sequence $Z_n=(z_k^n)_{k}\in\bb_\mathfrak{X}^\infty(Z)$. By
normalizing  we may suppose that for every $k,\ \|z^n_k\|=1$ and so
by our assumptions for $\widetilde{Z}$ and $Z$  there exists
$\widetilde{Z}_n=(\widetilde{z}_k^n)_{k}\in\bb_\mathfrak{D}^\infty(\widetilde{Z})$
such that $dist(Z_n,\widetilde{Z}_n)\leq\Delta/10C$. Then for all
$k$, $\|z_k^n-\widetilde{z}_k^n\|\leq \delta_k/10C$ and so
$\|\widetilde{z}_k^n\|\geq 1- \delta_k/10C$.  Let $k_0\geq n$ be such that
$x_{n-1}<z_{k_0}^n$ and let player I  play
$\widetilde{Z}_n|_{[k_0,\infty]}=(\widetilde{z}_k^n)_{k\geq k_0}$ in
the $n^{th}$- move of the discrete game
$G_\mathfrak{D}(\widetilde{Z})$. Then player II extends
$(\widetilde{x}_0,...,\widetilde{x}_{n-1})$ according to his
strategy in $G_\mathfrak{D}(\widetilde{Z})$ for $\mathcal{G}$, by
picking $\widetilde{x}_n\in<(\widetilde{z}_k^n)_{k\geq
k_0}>_\mathfrak{D}$. Then $\widetilde{x}_n=\sum_{k\in
I_n}\lambda_k^n\widetilde{z}_k^n$, where $I_n$ is a finite segment
in $\nn$ with $\min I_n\geq k_0$ and $\lambda_k^n\in\rr$. Going back
to Gowers'game $G_\mathfrak{X}(Z)$ let player II play
$x_n=\sum_{k\in I_n}\lambda_k^nz_k^n$. Then $x_n>x_{n-1}$ and so
player II    forms in this way a block sequence in
$\bb_\mathfrak{X}(Z)$.

 It remains to show that
$(x_n)_n\in\ff_\Delta$. Since $(\widetilde{x}_n)_n\in
\mathcal{G}\subseteq \ff_{\Delta/10C}\subseteq
(\bb^\infty_{B_\mathfrak{X}})_{\Delta/10C}$, we have that   for all $n$,
$\|\widetilde{x}_n\|\leq 1+\delta_n/10C$. Hence
\[|\lambda_k^n|\leq 2C\frac{\|\widetilde{x}_n\|}{\|\widetilde{z}_k^n\|}\leq
           2C \frac{1+\delta_n/10C}{1-\delta_k/10C}\leq 2C\frac{1+\delta_0/10C}
            {1-\delta_0/10C}\leq 4C,\]
 for all $k\in I_n$.

Therefore,  $\|x_n-\widetilde{x}_n\|\leq\sum_{k\in I_n}
            |\lambda_k^n|\|z_k^n-\widetilde{z}_k^n\|\leq 4C\sum_{k\in I_n}\frac{\delta_k}{10C}
            \leq \frac{4}{5}\delta_{\min I_n}\leq \frac{4}{5} \delta_n.$
Since  $(\widetilde{x}_n)_{n}\in\ff_{\Delta/10C}$, the last
inequality gives that $(x_n)_{n\in\nn}\in
\ff_{\frac{4\Delta}{5}+\frac{\Delta}{10C}}\subseteq\ff_\Delta$.\end{proof}

The above lemmas lead us to define the next property for  a subset $\mathfrak{D}$ of $ \mathfrak{X}$ and a given
 sequence $\Delta=(\delta_n)_n$ of positive real numbers.
\begin{enumerate}
\item[$(\mathfrak{D} 3)$] ($\Delta-$ \textit{block
covering property}) For every $\widetilde{Z}\in\bb_\mathfrak{D}^\infty$  there exists
$Z\in\bb_\mathfrak{X}^\infty$ such that $\bb_{B_\mathfrak{X}}^\infty
(Z)\subseteq(\bb_\mathfrak{D}^\infty(\widetilde{Z}))_\Delta$.
\end{enumerate}
In the next proposition we  give an example of a  subset
$\mathfrak{D}$ of $\mathfrak{X}$ with  properties
$(\mathfrak{D}1)-(\mathfrak{D}3)$. Actually we show that  a much
stronger than $(\mathfrak{D}3)$ property can be  satisfied. In
particular for every $\widetilde{Z}\in\bb_\mathfrak{D}^\infty$,
$\widetilde{Z}=(\widetilde{z}_n)_n$  setting $Z=(z_n)_n$ with
$z_n=\widetilde{z}_{2n}+\widetilde{z}_{2n+1}$ then
$\bb_{B_\mathfrak{X}}^\infty
(Z)\subseteq(\bb_\mathfrak{D}^\infty(\widetilde{Z}))_\Delta$.
\begin{prop} \label{prop16}For every  sequence $\Delta=(\delta_n)_{n}$ of positive real
numbers there is $\mathfrak{D}\subseteq B_\mathfrak{X}\cap <(e_n)_n>$ satisfying
$(\mathfrak{D}1)-(\mathfrak{D}3)$ and such that $(e_n)_n\in
\bb^\infty_\dd$.
\end{prop}
\begin{proof} Let $(k_n)_{n}$ be a strictly increasing sequence
of positive integers  such that for every $n$,
$2^{{-k_n}+1}\leq\delta_n$. For $i, l\in\nn$,  $l\geq 1$, let
\[\Lambda (i,l)=\{t\cdot 2^{-l\cdot(k_i+1)}:\;t\in\mathbb{Z}\}\]
For every finite nonempty segment $I=[n_1,n_2]$ of $\nn$, $n_1\leq
n_2$, define $\mathfrak{D}(I)=\mathfrak{D}([n_1,n_2])$ to be the set
of all $x=\sum_{i=n_1}^{n_2}\lambda_ie_i$ satisfying the following
properties.
 \begin{enumerate}
 \item[(i)] For all $n_1\leq i\leq n_2,$ $\lambda_i\in\Lambda (i,l)$,
where $l=n_2-n_1+1$ is the length of $I$.
 \item[(ii)] The coefficients $\lambda_{n_1}$ and $\lambda_{n_2}$ are both
 nonzero.
 \item[(iii)] $\|x\|\leq 1$.
\end{enumerate}
Finally we set \[\dd=\bigcup_{n_1\leq n_2}\dd([n_1,n_2])\]
 It is  easy to see that
$\mathfrak{D}$ satisfies $(\mathfrak{D}1)-(\mathfrak{D}2)$. In
particular $(e_n)_n\in \bb_\dd^\infty$. It remains to show that
$\mathfrak{D}$ has the $\Delta$- block covering property. Actually
we will
 prove that $\mathfrak{D}$ has a stronger property and to  do this
we first state the following.

\bigskip

\noindent\textit{Claim.} Let $\widetilde{Z}\in\bb_\dd^\infty$ and
let $w\in <\widetilde{Z}>$ such that
card$(\text{supp}_{\widetilde{Z}}(w))\geq 2$ and
$\|w\|\leq 1$.
 Then there
is $\widetilde{w}\in <\widetilde{Z}>_\dd$ such that
\begin{enumerate}
\item[(1)]$\text{supp}_{\widetilde{Z}}(\widetilde{w})=\text{supp}_{\widetilde{Z}}(w)$.
\item[(2)] $\|w-\widetilde{w}\|\leq 2^{-k_{m_1}+1}$, where $m_1=\min supp_{\widetilde{Z}}(w)$.
\end{enumerate}

\bigskip

\noindent\textit{Proof of the claim.} Let
$\widetilde{Z}=(\widetilde{z}_j)_j$ and let  $(I_j)_j$,
$I_j=[n_1(j),n_2(j)]$, $n_1(j)\leq n _2(j)$,  be the sequence of
successive finite nonempty segments of $\nn$ such that
$\widetilde{z}_j\in \dd(I_j)$. Let $m_1<m_2$ in $\nn$ and
$(\mu_j)_{j=m_1}^{m_2}$ be scalars such that $\mu_{m_1}$,
$\mu_{m_2}$ are both nonzero and let $w=\sum_{j\in[m_1, m_2]}\mu_j
\widetilde{z}_j$ in $B_\mathfrak{X}$.

Set $w'=(1- 2^{-k_{m_1}})w=\sum_{j\in[m_1, m_2]}(1-
2^{-k_{m_1}})\mu_j \widetilde{z}_j$ and
$\widetilde{w}=\sum_{j\in[m_1,m_2]}\widetilde{\mu}_j\widetilde{z}_j$,
where  $\widetilde{\mu}_j=s_j\cdot 2^{-(k_{n_1(j)}+1)}$ and if
$\mu_j\geq 0$, $ s_j=\lceil (1- 2^{-k_{m_1}})
\mu_j2^{k_{n_1(j)}+1}\rceil$ while  if $\mu_j<0,$ $ s_j=\lfloor (1-
2^{-k_{m_1}}) \mu_j2^{k_{n_1(j)}+1}\rfloor$, i.e. $\widetilde{\mu}_j$
are of the form $s_j\cdot2^{-(k_{n_1(j)+1})}$ such that $|\widetilde{\mu}_j|\geq|\mu_j(1-2^{-k_{m_1}})|$
and $|\widetilde{\mu}_j-(1-2^{-k_{m_1}})\mu_j|<2^{-(k_{n_1(j)}+1)}$.

It is easy to see that  $\widetilde{\mu}_j=0$ if and only if
$\mu_j=0$ and so
$\text{supp}_{\widetilde{Z}}(\widetilde{w})=\text{supp}_{\widetilde{Z}}
(w)$. Moreover for all $j$, $
|(1-2^{-k_{m_1}})\mu_j-\widetilde{\mu}_j|\leq 2^{-(k_{n_1(j)}+1)}$
and so
\begin{equation}\label{eq1}\begin{split}\|w'-\widetilde{w}\|&\leq\sum_{j\in[m_1,m_2]}\big|(1-2^{-k_{m_1}})\mu_j-\widetilde{\mu}_j\big|
\|\widetilde{z}_j\|\\&\leq\sum_{j\in[m_1,m_2]}2^{-(k_{n_1(j)}+1)}\leq
2^{-k_{n_1(m_1)}}\end{split}\end{equation} and therefore, since
$m_1\leq n_1(m_1) $, $\|w'-\widetilde{w}\|\leq 2^{-k_{m_1}}$.  As
$\|w-w'\|\leq 2^{-k_{m_1}}$, we obtain that $\|w
-\widetilde{w}\|\leq 2^{-k_{m_1}+1}$.
\\It remains to show that $\widetilde{w}\in\mathfrak{D}$. Since  for
all $j\in [m_1,m_2]$, $\widetilde{z}_j\in\dd(I_j)$,
 we have that $\widetilde{z}_j=\sum_{i\in
I_j}t^j_i 2^{-l_j(k_{i}+1)}e_i$, where  $l_j=n _2(j)-n _1(j)+1$ is
the length of $I_j$ and $t^j_{n_1(j)},t^j_{n_2(j)}$ are both
nonzero. Therefore setting $I=[n_1(m_1),n_2(m_2)]$, we have that
\begin{equation}\label{eq2}\widetilde{w}=\sum_{j\in[m_1,m_2]}
\widetilde{\mu}_j\widetilde{z}_j=\sum_{j\in[m_1,m_2]}\widetilde{\mu}_j\big(\sum_{i\in
I_j}t^j_i 2^{-l_j(k_{i}+1)}e_i\big)=\sum_{i\in
I}\lambda_ie_i\end{equation} where for all $i\in I_j$ and  $j\in
[m_1,m_2]$, $\lambda_i=t^j_i 2^{-l_j(k_{i}+1)}\widetilde{\mu}_j$ and
 $\lambda_i=0$, for all $i\in I\setminus \bigcup_{j\in[m_1,m_2]}I_j$.

We first show that  condition (i) of the definition of
$\mathfrak{D}$ is satisfied, that is  for all $i\in I$,
$\lambda_i\in\Lambda(i,l)$ where $l=n_2(m_2)-n_1(m_1)+1$ is the
length of $I$. Since $0\in \Lambda(i,l)$,  it suffices to check it
for each $i\in\bigcup_{j\in[m_1,m_2]}I_j$. So fix $j\in[m_1,m_2]$
and $i\in I_j$. Then
 \begin{equation}\lambda_i=t^j_i 2^{-l_j(k_{i}+1)}\widetilde{\mu}_j=t^j_i 2^{-l_j(k_{i}+1)}  s_j
2^{-(k_{n_1(j)}+1)}=\tau^j_i{2^{-l(k_i+1)}}\end{equation} where
$\tau^j_i=t^j_i s_j 2^{(l-l_j)(k_i+1)-(k_{n_1(j)}+1)}$. Since
$m_1<m_2$ we have that $l> l_j$. Also $n_1(j)\leq i$ and so
$(l-l_j)(k_i+1)-(k_{n_1(j)}+1)\geq 0$. Therefore  $\tau^j_i\in
\mathbb{Z}$ which gives that $\lambda_i\in \Lambda(i,l)$.

Moreover, since $\widetilde{\mu}_{m_1},\widetilde{\mu}_{m_2},
t_{n_1(m_1)}^{m_1}, t_{n_2(m_2)}^{m_2}$ are all non zero we have
that $\lambda_{n_1(m_1)}$ and $\lambda_{n_2(m_2)}$ are also non zero
and so condition (ii) of the definition of $\mathfrak{D}$ is also
satisfied. Finally by (\ref{eq1}),  $\|\widetilde{w}\|\leq
\|w'\|+2^{-k_{n_1(m_1)}}\leq 1$ and so condition (iii) is fulfilled.
By the above we have that  $\widetilde{w}\in \dd$ and the proof of
the claim is complete.

We continue with the proof of the proposition. Let
$\widetilde{Z}=(\widetilde{z}_j)_j$ in $\bb_\dd^\infty$ and let
$Z=(z_j)_j$ where for all $j$,
$z_j=\widetilde{z}_{2j}+\widetilde{z}_{2j+1}$. Pick $W=(w_i)_i$ in
$\bb_{B_\mathfrak{X}}^\infty(Z)$. Then for each $i$ there exist
$m_1^i<m_2^i$ and scalars $(\mu_j)_j$ such that
$w_i=\sum_{j\in[m_1^i, m_2^i]}\mu_j \widetilde{z}_j\in
B_\mathfrak{X}$ and $\mu_{m_1^i},\mu_{m_2^i}$ are both non zero. By
the claim, for each $i$ there exist scalars $(\widetilde{\mu}_j)_j$
such that $\widetilde{w}_i=\sum_{j\in[m_1^i,
m_2^i]}\widetilde{\mu}_j \widetilde{z}_j\in \dd$ and
$\|w_i-\widetilde{w}_i\|\leq 2^{{-k_{m_1^i}}+1}\leq 2^{{-k_i}+1}\leq
\delta_i$. We set $\widetilde{W}=(\widetilde{w}_i)_i$ and then
$\widetilde{W}\in \bb_\dd^\infty(\widetilde{Z})$ and
$dist(\widetilde{W},W)\leq \Delta$. Hence
$\bb_{B_\mathfrak{X}}^\infty(Z)\subseteq
(\bb^\infty_\dd(\widetilde{Z}))_\Delta$ and the proof is complete.
\end{proof}

It is easy to see that
$\rho(x,y)=\|x-y\|+|\frac{1}{\|x\|}-\frac{1}{\|y\|}|$, $x,y\in
\mathfrak{X}\setminus\{0\}$ is an equivalent metric on
$(\mathfrak{X}\setminus\{0\},\|\cdot\|)$ and that  the product
topology on $(\mathfrak{X}\setminus\{0\},\rho)^\nn$ makes
$\bb_{\mathfrak{X}}^\infty$
  a Polish space.

\begin{lem} \label{analytikotita} Let $\ff$ be an  analytic subset
of $\bb_{\mathfrak{X}}^\infty$ and  $\Delta=(\delta_n)_{n}$
be a sequence of positive real numbers. Then
\begin{enumerate}
\item[(i)] $\ff_\Delta$ is  analytic in
$\bb_{\mathfrak{X}}^\infty$.
\item[(ii)] For every countable $\mathfrak{D}\subseteq \mathfrak{X}$, $\ff_\Delta\cap\bb_\mathfrak{D}^\infty$ is analytic in
$\mathfrak{D}^\nn$ (where $\mathfrak{D}$ is endowed with the
discrete topology).
\end{enumerate}
\end{lem}
\begin{proof}
\noindent (i) It is easy to see that   $\mathcal{Q}=\{(U,V):
dist(U,V)\leq\Delta\}$ is closed in
$\bb_{\mathfrak{X}}^\infty\times \bb_{\mathfrak{X}}^\infty$. Let
$proj_1$  (resp. $proj_2$) be the projection of
$\bb_{\mathfrak{X}}^\infty\times \bb_{\mathfrak{X}}^\infty$ onto
the first (resp. second) coordinate. Then notice that
$\ff_\Delta=proj_1[\mathcal{Q}\cap(\bb_{\mathfrak{X}}\times
\ff)]=proj_1[\mathcal{Q}\cap proj_2^{-1}(\ff)]$.

\noindent (ii) Let $I:\mathfrak{D}^\nn\to \mathfrak{X}^\nn$
 be the identity
map. Then $I$ is clearly continuous   and
 $\ff_\Delta\cap\bb_\mathfrak{D}^\infty=I^{-1}(\ff_\Delta)$.
 \end{proof}
\begin{thm}(W. T. Gowers)\label{main_Gowers_combinatorial_theorem}
Let $\mathfrak{X}$ be a normed linear space with a basis and let
$\ff\subseteq \bb_{B_\mathfrak{X}}^\infty$ be an analytic family of
block sequences in the unit ball  $B_\mathfrak{X}$ of
$\mathfrak{X}$. Then for every  $\Delta>0$ there exists a block
sequence $Z\in \bb_\mathfrak{X}^\infty$ such that either
$\bb_{B_\mathfrak{X}}^\infty(Z)\cap\ff=\emptyset$ or player II has
 a winning strategy in Gowers' game $G_{\mathfrak{X}}(Z)$ for
$\ff_\Delta$.
\end{thm}
\begin{proof} Let $(e_n)_n$ be a normalized basis for $\mathfrak{X}$ with constant $C$.
Let  $\Delta'=(\delta'_n)_n$ be a sequence of positive real numbers such that  $\delta'_0\leq 1$,  $\delta'_n\leq \delta_n$,
 and
$\sum_{i>n}\delta'_i\leq \delta'_n$.
 By Proposition \ref{prop16}, there is $\mathfrak{D}\subseteq
\mathfrak{X}$ with $(e_n)_n\in\bb_\dd^\infty$  satisfying
$(\mathfrak{D}1)-(\mathfrak{D}3)$  for $\Delta'/10C$. Let also
$\mathcal{G}=\ff_{\Delta'/10C }\cap\bb_\mathfrak{D}^\infty$. By
Lemma \ref{analytikotita}, $\mathcal{G}$ is analytic in
$\mathfrak{D}^\nn$ and applying  Theorem
\ref{main_cobinatorial_theorem}, we obtain a block sequence
$\widetilde{Z}\in\bb_\dd^\infty$ such that  either
$\bb_\dd^\infty(\widetilde{Z})\cap \mathcal{G}=\emptyset$ or player
II has winning strategy in $G_\dd(\widetilde{Z})$ for $\mathcal{G}$.
Choose $Z\in\bb_{\mathfrak{X}}^\infty$ such that
$\bb_{B_\mathfrak{X}}^\infty(Z)\subseteq
(\bb^\infty_\dd(\widetilde{Z}))_{\Delta'/10C}$. From Lemmas
\ref{lemma2_from_discrete_to_continuous} and
\ref{lemma1_from_discrete_to_continuous}, we have that either
$\bb_{B_\mathfrak{X}}^\infty(Z)\cap
\mathcal{F}=\emptyset$, or player II has a winning strategy in Gowers'
game $G_{\mathfrak{X}}(Z)$ for $\mathcal{F}_{\Delta'}$ and so  (as $\Delta'\leq \Delta$) for
$\mathcal{F}_{\Delta}$ as well.
\end{proof}

\section{A Ramsey consequence on  $k$-tuples of block bases.}
The main goal of  this section is to  prove  Theorem \ref{Ramsey}.
First we need to do  some preliminary work and introduce some
notation . Fix a positive integer $k\geq 2$. For each $0\leq i\leq
k-1$ and every infinite subset $L=\{l_0<l_1<...\}$ of $\nn$ we set
$L_{i(mod k)}=\{l_{k n+i}:n\in\nn\}$ and we define
\[([L]^\infty)^k_\circ=\prod_{i=0}^{k-1}[L_{i(mod k)}]^\infty=\{(L_i)_{i=0}^{k-1}\in([L]^\infty)^k:\forall i\; L_i\subseteq L_{i(mod k)}\}\]
Notice that  $([L]^\infty)_\circ^k$ is not hereditary, that is
generally  $([L']^\infty)^k_\circ\nsubseteq([L]^\infty)^k_\circ$,
for $L'\subseteq L$. Let also
 \[([L]^\infty)^k_\perp=\{(L_i)_{i=0}^{k-1}\in([L]^\infty)^k:\forall i\neq j\; L_i\cap L_j=\emptyset\}\]
We have the following elementary lemma which relates the above types
of products.
\begin{lem}
 \label{h11} Let   $N=\{(2n+1)k:n\in\nn\}$. Then
$([N]^\infty)_\perp^k\subseteq
\bigcup_{L\in[\nn]^\infty}([L]^\infty)_{\circ}^k$.
\end{lem}
\begin{proof} Let  $(M_i)_{i=0}^{k-1}\in([N]^\infty)^k_\perp$. Let
$M=\bigcup_{i=0}^{k-1} M_i$  and for each $m\in M$  define the
interval $I_m=[m-i_m,m-i_m+k-1]$  of $\nn$  where $i_m$ is the
unique natural number $i$ such that $m\in M_i$. Notice that the
length of all $I_m$ is $k$ while the length of an interval with
nonequal endpoints in $N$ is at least $2k+1$. Hence  for  $m_1\neq
m_2$, $I_{m_1}\cap I_{m_2}=\emptyset$ and  for all $m\in M$,
$I_m\cap N=\{m\}$.

Let $L=\bigcup_{m\in M}I_m$. We claim that $(M_i)_{i=0}^{k-1}\in
([L]^\infty)_\circ^k$. Indeed, let $L=(l_n)_n$ be the increasing
enumeration of $L$. For each $0\leq i\leq k-1$ and $m\in M$ let
 $I_m(i)=m-i_m+i$ be the $i^{th}$-element
of $I_m$. Since $(I_m)_{m\in M}$ is a sequence of pairwise disjoint
intervals of $\nn$ of length $k$, we easily see that $L_{i(mod
k)}=\bigcup_{m\in M}I_m(i)$. Fix $0\leq i\leq k-1$. Then $m\in M_i$
if and only if $i_m=i$ if and only if $I_m(i)=m$. Hence
$M_i=\bigcup_{m\in M_i}\{I_m(i)\}\subseteq \bigcup_{m\in
M}\{I_m(i)\}= L_{i(mod k)}$.\end{proof}

The above notation is easily extended to block sequences in the
unit ball $B_\mathfrak{X}$ of a Banach space $\mathfrak{X}$ as
follows. For every $Z\in \bb_\mathfrak{X}^\infty$ let
\[(\bb_{B_\mathfrak{X}}^\infty(Z))_\circ^k=\{(Z_i)_{i=0}^{k-1}\in
(\bb_{B_\mathfrak{X}}^\infty)^k:\;\forall i \;\;Z_i\preceq
Z|_{\nn_{i(mod k)}}\}\] and generally for $L\in [\nn]^\infty$, we
set
\[(\bb_{B_\mathfrak{X}}^\infty(Z|_L))_\circ^k=\{(Z_i)_{i=0}^{k-1}\in
(\bb_{B_\mathfrak{X}}^\infty)^k:\;\forall i \;\;Z_i\preceq
Z|_{L_{i(mod k)}}\}\] The next lemma is an immediate consequence of Lemma \ref{h11}.
\begin{lem}
 \label{h1} Let $Z\in \bb_\mathfrak{X}^\infty$ and $N=\{(2n+1)k:n\in\nn\}$. Then
\[(\bb_{B_\mathfrak{X}}^\infty(Z|_N))_\perp^k\subseteq
\bigcup_{L\in[\nn]^\infty}(\bb_{B_\mathfrak{X}}^\infty(Z|_L))_{\circ}^k.\]\end{lem}

For a family  $\mathfrak{F}\subseteq
(\bb_{B_\mathfrak{X}}^\infty)^k$ let
\[\mathcal{F}^\mathfrak{F}=\{Z\in
\bb_{S_\mathfrak{X}}^\infty:\;\mathfrak{F}\cap
(\bb_{B_\mathfrak{X}}^\infty(Z))_{\circ}^k\neq \emptyset\},\]
where $S_\mathfrak{X}$ is the unit sphere of $\mathfrak{X}$.

\begin{lem}\label{h2}  If  $\mathfrak{F}$ is   analytic in $(\bb_{\mathfrak{X}}^\infty)^k$,
then $\ff^\mathfrak{F}\subseteq \bb_{S_\mathfrak{X}}^\infty$ is
analytic in $\bb_{\mathfrak{X}}^\infty$.\end{lem}
\begin{proof} Let $\mathcal{K}=\{(Z,(V_i)_{i=0}^{k-1})\in \bb_{S_\mathfrak{X}}^\infty \times
(\bb_{B_\mathfrak{X}}^\infty)^k:\;(V_i)_{i=0}^{k-1}\in
(\bb_{B_\mathfrak{X}}^\infty(Z))_\circ^k\}$.  Then $\mathcal{K}$
is a closed subset of $\bb_{\mathfrak{X}}^\infty \times
(\bb^\infty_{{\mathfrak{X}}})^k$ and that
$\mathcal{F}^\mathfrak{F}=proj_1\big[\big(\bb_{\mathfrak{X}}^\infty\times\mathfrak{F}\big)\cap
\mathcal{K}\big]$. \end{proof}

\begin{proof}[\textbf{Proof of Theorem \ref{Ramsey}: }]Let $(e_n)_n$ be a
normalized basis of $\mathfrak{X}$ with basis constant $C$. Choose
$\Delta'=(\delta'_n)_n$ such that $0<\delta'_n\leq (4
C)^{-1}\delta_n$ and $\ \sum_{j=n+1}^\infty\delta'_j\leq
\delta'_n$. By Lemma \ref{h2}, we have that
$\mathcal{F}^\mathfrak{F}$ is an analytic subset of
$\bb_{B_\mathfrak{X}}^\infty$ and  by Theorem
\ref{main_Gowers_combinatorial_theorem} there is a block
subsequence $Z=(z_n)_n$ such that either
$\bb_{B_\mathfrak{X}}^\infty(Z)\cap\ff^\mathfrak{F}=\emptyset$ or
player II has winning strategy in Gowers' game
$G_{\mathfrak{X}}(Z)$ for $(\ff^\mathfrak{F})_{\Delta'}$. Let
$Y=Z|_N$, where $N=\{(2n+1)k:n\in\nn\}$. We claim that $Y$
satisfies the conclusion of the theorem.

Indeed, if
$\bb_{B_\mathfrak{X}}^\infty(Z)\cap\ff^\mathfrak{F}=\emptyset$ then
for all $Z'\in \bb_{B_\mathfrak{X}}^\infty(Z)$,
$\mathfrak{F}\cap(\bb_{B_\mathfrak{X}}^\infty(Z'))_{\circ}^k=\emptyset$.
In particular for all $L\in[\nn]^\infty$,
$\mathfrak{F}\cap(\bb_{B_\mathfrak{X}}^\infty(Z|_L))_{\circ}^k=\emptyset$
which by Lemma \ref{h1} gives that
$\mathfrak{F}\cap(\bb_{B_\mathfrak{X}}^\infty(Y))_\perp^k=\emptyset$.

So let us assume  that   player II has a winning strategy in Gowers'
game $G_{\mathfrak{X}}(Z)$ for $(\ff^\mathfrak{F})_{\Delta'}$. Since $Y=Z|_N$  the
same holds for the game $G_{\mathfrak{X}}(Y)$. Fix
$(U_i)_{i=0}^{k-1}\in(\bb_{B_\mathfrak{X}}^\infty(Y))^k$. We have to
show that there exists $(V_i)_{i=0}^{k-1}\in(
\bb^\infty_\mathfrak{X})^k$ such that  $V_i\preceq U_i$ and
$(V_i)_{i=0}^{k-1}\in \mathfrak{F}_\Delta$. Consider a run of the
game such that in the $n^{th}$- move player I plays $U_i$, where
$n=i(modk)$. Then player II succeeds to construct a block sequence
$V=(v_n)_n$ in $(\mathcal{F}^\mathfrak{F})_{\Delta'}$ such that
$v_n\in U_i$ for all $n=i(modk)$.  Choose $W$ in
$\mathcal{F}^\mathfrak{F}$ with $dist(V, W)\leq \Delta'$ and for
each $i$,  $W_i\preceq W|_{\nn_{i(modk)}}$ such that
$(W_i)_{i=0}^{k-1}\in(\bb_{B_\mathfrak{X}}^\infty(W))^k_{\circ}\cap\mathfrak{F}$.
Let $W=(w_n)_n$ and $W_i=(w^i_n)_n$. Then for each $i=1,...,k$ there is a block sequence
$(F_n^i)_n$ of finite subsets of $\nn_{i(modk)}$ and a sequence of
scalars $(\lambda_j)_j$ such that for all $i$ and all $n$,
$w^i_n=\sum_{j\in F^i_n}\lambda_j w_j$. We set $v^i_n=\sum_{j\in
F^i_n}\lambda_j v_j$ and let $V_i=(v^i_n)_n$. Then for all $i$,
$V_i\preceq V|_{\nn_{i(mod k)}}\preceq U_i$. It remains to show that
$(V_i)_{i=0}^{k-1}\in \mathfrak{F}_\Delta$. For this it suffices
to see that $dist(V_i,W_i)\leq \Delta$, for all $i$. Indeed fix
$0\leq i\leq k-1$ and $n\in\nn$. Since $\|w^i_n\|\leq 1$ and
$\|w_j\|=1$ , we get that $|\lambda_j|\leq 2C$ and therefore
\[\|v_n^i-w^i_n\|\leq
\sum_{j\in F^i_n}|\lambda_j|\|v_j-w_j\|\leq 2C \sum_{j\in
F^i_n}\delta'_j\leq 4C\delta'_n\leq \delta_n\] Hence
$(U_i)_{i=0}^{k-1}\in(\mathfrak{F}_\Delta)^\uparrow$.
\end{proof}

\section{Comments}

\noindent \textbf{1.} C. Rosendal in \cite{R1} proves a Ramsey
dichotomy between  winning strategies in Gowers' game and winning
strategies in the infinite asymptotic game.  By appropriately
modifying his argument, one can check that the proof in \cite{R1}
works in the more general setting of a linear space $\mathfrak{X}$
of countable dimension over the field of reals  provided that both
games are restricted on a \textit{countable} subset $\mathfrak{D}$
of $\mathfrak{X}$ satisfying property $(\mathfrak{D}1)$ stated in
the introduction. This modification can be used to derive an
alternative proof of Theorem \ref{main_cobinatorial_theorem}.
\medskip

\noindent \textbf{2.}  Theorem \ref{Ramsey} is actually an extension of the
following fact concerning pairs of infinite subsets of $\nn$. Given
an analytic family $\mathfrak{F}\subseteq [\nn]^\infty\times
[\nn]^\infty$ there is an infinite subset $L$ of $\nn$ such that
either all \textit{disjoint} pairs of infinite subsets of $L$ belong
to the complement of $\mathfrak{F}$ or for every $(L_1,L_2)\in
[L]^\infty\times [L]^\infty$, there is $(L'_1,L'_2)\in \mathfrak{F}$
such that $L'_i\subseteq L_i$ for all $i=1,2$ . To see this consider
the map $\Phi: M\to (M_0,M_1)$ where if $M=\{m_i\}_i$ is the
increasing enumeration of $L$ then $M_0=\{m_i\}_{i\; \text{even}}$
and $M_1=\{m_i\}_{i \;\text{odd}}$. Then apply Silver's theorem (see
\cite{S}) for the family $\Phi^{-1}(\mathfrak{F}^\uparrow)$ where
$\mathfrak{F}^\uparrow=\{(L,M): \exists (L',M')\in\mathfrak{F}\;\;
\text{with}\; L'\subseteq L\;\text{and}\; M'\subseteq M\}$. It is
easy to see that keeping the ``half" of the monochromatic set the
result follows. Also, applying K. Milliken's theorem \cite{Mil},
one can derive an analogue of the above result for pairs of block
sequences of finite subsets of $\nn$.

\section{Acknowledgement}
We would like to thank the referee for his (or her) suggestions
which simplified the proof of Theorem
\ref{main_cobinatorial_theorem}.

\end{document}